\documentclass{conm-p-l}

\newtheorem{theorem}[equation]{Theorem}
\newtheorem{lemma}[equation]{Lemma}
\newtheorem{corollary}[equation]{Corollary}
\newtheorem{proposition}[equation]{Proposition}

\newcommand{\updot}{\textstyle\cdot}

\numberwithin{equation}{section}

\begin{document}

\title{Exponential sums on ${\bf A}^n$, III}

% author one information
\author{Alan Adolphson}
\address{Department of Mathematics, Oklahoma State University,
  Stillwater, Oklahoma 74078}
\email{adolphs@math.okstate.edu}
\thanks{The first author was supported in part by NSA Grant
  \#MDA904-97-1-0068}

% author two information
\author{Steven Sperber}
\address{School of Mathematics, University of Minnesota, Minneapolis,
  Minnesota 55455}
\email{sperber@math.umn.edu}

\date{\today}

\begin{abstract}
  We give two applications of our earlier work\cite{AS3}.  We compute
  the $p$-adic cohomology of certain exponential sums on ${\bf A}^n$
  involving a polynomial whose homogeneous component of highest degree
  defines a projective hypersurface with at worst weighted homogeneous
  isolated singularities.  This study was motivated by recent work of
  Garc\'{\i}a\cite{G}.  We also compute the $p$-adic cohomology of
  certain exponential sums on ${\bf A}^n$ whose degree is divisible by
  the characteristic.
\end{abstract}

\maketitle

\section{Introduction}

Let $p$ be a prime number, $q=p^a$, and ${\bf F}_q$ the finite field
of $q$ elements.  Associated to a polynomial $f\in{\bf
  F}_q[x_1,\ldots,x_n]$ and a nontrivial additive character $\Psi:{\bf
  F}_q\rightarrow{\bf C}^{\times}$ are exponential sums
\begin{equation}
S({\bf A}^n({\bf F}_{q^i}),f)=\sum_{x_1,\ldots,x_n\in{\bf F}_{q^i}}
\Psi({\rm Trace}_{{\bf F}_{q^i}/{\bf F}_q}f(x_1,\ldots,x_n))
\end{equation}
and an $L$-function
\begin{equation} L({\bf A}^n,f;t)=\exp\biggl(\sum_{i=1}^{\infty}
  S({\bf A}^n({\bf F}_{q^i}),f)\frac{t^i}{i}\biggr).
\end{equation}
Dwork has associated to $f$ a complex
$(\Omega^{\textstyle\cdot}_{C(b)},D)$ (of length $n$), depending on a
choice of rational parameter $b$ satisfying $0<b<p/(p-1)$ (see
\cite{AS4} for details).  Each cohomology group
$H^i(\Omega^{\textstyle\cdot}_{C(b)},D)$ is a vector space over a
field $\tilde{\Omega}_0$ (a finite extension of ${\bf Q}_p$) and has a
Frobenius operator $F$ satisfying 
\[ L({\bf A}^n,f;t)=\prod_{i=0}^n \det(I-tF\mid
H^i(\Omega^{\textstyle\cdot}_{C(b)},D))^{(-1)^{i+1}}. \] 

We write ${\bf F}_q[x]$ for ${\bf F}_q[x_1,\ldots,x_n]$ and consider
the complex $(\Omega^{\updot}_{{\bf F}_q[x]/{\bf F}_q},\phi_f)$, 
where $\Omega^k_{{\bf F}_q[x]/{\bf F}_q}$ denotes the module of
differential $k$-forms of ${\bf F}_q[x_1,\ldots,x_n]$ over ${\bf F}_q$
and $\phi_f:\Omega^k_{{\bf F}_q[x]/{\bf F}_q}\rightarrow\Omega^{k+1}_{{\bf
    F}_q[x]/{\bf F}_q}$ is defined by
\[ \phi_f(\omega)=df\wedge\omega, \]
where $d:\Omega^k_{{\bf F}_q[x]/{\bf
    F}_q}\rightarrow\Omega^{k+1}_{{\bf F}_q[x]/{\bf F}_q}$ is the
exterior derivative.  Every $\omega\in\Omega^k_{{\bf
  F}_q[x]/{\bf F}_q}$ can be uniquely written in the form 
\[ \omega=\sum_{1\leq i_1<\cdots<i_k\leq n} \omega(i_1,\ldots,i_k)\,
dx_{i_1}\wedge\cdots\wedge dx_{i_k}, \]
with $\omega(i_1,\ldots,i_k)\in{\bf F}_q[x]$.  If each coefficient
$\omega(i_1,\ldots,i_k)$ is a homogeneous form of degree $l$, we call
$\omega$ {\it homogeneous\/} and define
\[ \deg\omega =l+(n-k)(\delta-1), \]
where $\delta=\deg f$.  The point of this definition is that we can
define an increasing filtration $F.$ on
$\Omega^k_{{\bf F}_q[x]/{\bf F}_q}$ by setting
\[ F_l\Omega^k_{{\bf F}_q[x]/{\bf F}_q} = \text{the ${\bf
  F}_q$-span of homogeneous $k$-forms $\omega$ with $\deg\omega\leq l$}, \]
and $(\Omega^{\updot}_{{\bf F}_q[x]/{\bf
    F}_q},\phi_f)$ then becomes a filtered complex.  Consider the 
associated spectral sequence 
\begin{equation}
E_1^{r,s}=H^{r+s}(F_r/F_{r-1}(\Omega^{\updot}_{{\bf
    F}_q[x]/{\bf F}_q},\phi_f))\Rightarrow
    H^{r+s}(\Omega^{\updot}_{{\bf F}_q[x]/{\bf F}_q},\phi_f).
\end{equation}
As an immediate consequence of \cite[Theorem 1.13]{AS3},
we have the following.
\begin{theorem}
  Suppose there exists a positive integer $e$ satisfying 
\begin{equation}
\biggl(1+\frac{p}{(p-1)^2}\biggr)(e-1)<\delta 
\end{equation}
such that $E^{r,s}_e=0$ for all $r,s$ with $r+s\neq n$.  Then for
\begin{equation}
\frac{\delta}{(p-1)(\delta-e+1)}<b<\frac{p\delta}{(p-1)\delta+e-1}
\end{equation}
we have
\begin{equation}
H^i(\Omega^{\textstyle\cdot}_{C(b)},D) =0 \quad\text{for $i\neq n$}
\end{equation}
and
\begin{equation}
\dim_{\tilde{\Omega}_0}H^n(\Omega^{\textstyle\cdot}_{C(b)},D) =M_f,
\end{equation}
where $M_f$ is the sum of the Milnor numbers of the critical points of
the mapping $f:{\bf A}^n\rightarrow{\bf A}^1$.  In particular,
$L({\bf A}^n,f;t)^{(-1)^{n+1}}$ is a polynomial of degree $M_f$.
\end{theorem}

{\it Remark}.  Inequality (1.5) is equivalent to the assertion that
the right-most term in (1.6) is greater than the left-most term in
(1.6), i.~e., (1.5) is equivalent to the existence of rational $b$
satisfying (1.6).  It is explained in \cite[section 1]{AS3} that the
vanishing of $E_e^{r,s}$ for all $r,s$ with $r+s\neq n$ implies that
$f:{\bf A}^n\rightarrow{\bf A}^1$ has isolated critical points, hence
the sum of the Milnor numbers is finite.  

To apply Theorem 1.4, one must find conditions on the polynomial $f$
that guarantee that, for some $e\geq 1$, $E_e^{r,s}=0$ for all $r,s$
with $r+s\neq n$.  When $e=1$, this is equivalent to the condition
that the partial derivatives of the homogeneous component of degree
$\delta$ of $f$ form a regular sequence in ${\bf F}_q[x]$.  When
$e>1$, the problem is much harder.  We gave one example of such a
condition in \cite[section 5]{AS3}.  The purpose of this article is to
give two more examples of such conditions.

Write
\begin{equation}
f=f^{(\delta)}+f^{(\delta')}+f^{(\delta'-1)}+\cdots+f^{(0)}, 
\end{equation}
where $f^{(i)}$ is homogeneous of degree $i$ and $1\leq \delta'\leq
\delta-1$, i.~e., $f^{(\delta')}$ is the homogeneous component of
second-highest degree of $f$.  We prove the following
result, which was stated in \cite{AS3}.  (The terms ``weighted
homogeneous'' isolated singularity and ``total degree'' of a weighted
homogeneous isolated singularity will be defined in the next section.)
\begin{theorem}
Suppose that the hypersurface $f^{(\delta)}=0$ in ${\bf P}^{n-1}$ has at
worst weighted homogeneous isolated singularities, of total degrees
$\delta_1,\ldots,\delta_s$, and that none of these singularities lies
on the hypersurface $f^{(\delta')}=0$ in ${\bf P}^{n-1}$.  Suppose also
that $(p,\delta\delta'\delta_1\cdots\delta_s)=1$.  Then
$E^{r,s}_{\delta-\delta'+1}=0$ for 
all $r,s$ with $r+s\neq n$.
\end{theorem}

The hypothesis of Theorem 1.10, with $\delta'=\delta-1$, was first
considered by Garc\'{\i}a\cite{G}.  He showed in that case that the $l$-adic
cohomology groups of the exponential sum (1.1) vanish except in degree
$n$ and that the reciprocal roots of $L({\bf A}^n,f;t)^{(-1)^{n+1}}$
are pure of weight $n$.  In particular, he obtains the estimate
\[ |S({\bf A}^n({\bf F}_{q^i}),f)|\leq M_f q^{ni/2}. \]
By our approach, we have not been able to obtain archimedian estimates
for the reciprocal roots of $L({\bf A}^n,f;t)^{(-1)^{n+1}}$.  Thus that
question is still open for $\delta'<\delta-1$, although we conjecture
the roots are again pure of weight $n$ when (1.5) holds for
$e=\delta-\delta'+1$. 

For our second example, we consider the case where the degree of $f$
is divisible by $p$.  This uses ideas similar to those in the proof of
Theorem 1.10, but the computations are simpler.  As we noted in
\cite{AS4}, if $\{\partial f^{(\delta)}/\partial x_i\}_{i=1}^n$ form a
regular sequence in ${\bf F}_q[x]$, then $L({\bf
  A}^n,f;t)^{(-1)^{n+1}}$ is a polynomial of degree $(\delta-1)^n$ all of
whose reciprocal roots have absolute value $q^{n/2}$, even if $p|\delta$.
We consider the case where these partial derivatives do not form a
regular sequence.
\begin{theorem}
  Suppose $p|\delta$ and the set of common zeros of $\{\partial
  f^{(\delta)}/\partial x_i\}_{i=1}^n$ in ${\bf P}^{n-1}$ is finite
  and nonempty.  Suppose also that $(p,\delta')=1$ and that the
  hypersurface $f^{(\delta')}=0$ in ${\bf P}^{n-1}$ contains no common
  zero of $\{\partial f^{(\delta)}/\partial x_i\}_{i=1}^n$.  Then
  $E^{r,s}_{\delta-\delta'+1}=0$ for all $r,s$ with $r+s\neq n$.
\end{theorem}

{\it Remark}.  For example, if $p|\delta$ and $f^{(\delta)}=0$ defines
a smooth hypersurface in~${\bf P}^{n-1}$, then the set of common
zeroes of $\{\partial f^{(\delta)}/\partial x_i\}_{i=1}^n$ in ${\bf
  P}^{n-1}$ is finite. (If the set of common zeros had dimension $\geq
1$, it would have nonempty intersection with the hypersurface
$f^{(\delta)}=0$, and any such point of intersection would be a
singular point of this hypersurface.) We conjecture that under the
hypothesis of Theorem~1.11, the reciprocal roots of $L({\bf
  A}^n,f;t)^{(-1)^{n+1}}$ are pure of weight $n$ when (1.5) holds for
$e=\delta-\delta'+1$. 

Garc\'{\i}a also gave a formula for the degree of $L({\bf
  A}^n,f;t)^{(-1)^{n+1}}$ in terms of the Milnor numbers of the
singularities of $f^{(\delta)}=0$.  We derive an analogous formula in
section 6 under the hypothesis of either Theorem 1.10 or~1.11.
  
For certain of the constructions made in the proofs of
Theorems~1.10 and~1.11, it may be necessary to extend scalars from ${\bf
  F}_q$ to a larger finite field.  Since such extensions of scalars do
not affect the computation of cohomology, we make no further comment
on them.
\section{Hypersurface singularities}

For this general discussion of singularities, we work over an
arbitrary algebraically closed field $K$.  
Let $f\in K[x_1,\ldots,x_n]$, put ${\bf
  0}=(0,\ldots,0)$, and assume $f({\bf 0})=0$.  We say that the
hypersurface $f=0$ has an {\it isolated singularity at\/} ${\bf 0}$ if
${\bf 0}$ is an isolated critical point of the map $f:K^n\rightarrow
K$, i.~e., there exists a Zariski open neighborhood $U$ of~${\bf 0}$
in~$K^n$ such that the only common zero of $\partial f/\partial
x_1,\ldots,\partial f/\partial x_n$ on $U$ is ${\bf 0}$.  Let ${\bf
  m}=(x_1,\ldots,x_n)$, the maximal ideal of $K[x_1,\ldots,x_n]$
corresponding to ${\bf 0}$ and let $K[x_1,\ldots,x_n]_{\bf m}$ be the
localization of $K[x_1,\ldots,x_n]$ at ${\bf m}$.  When ${\bf 0}$ is
an isolated singularity, then 
\begin{equation}
\mbox{Krull dim}\;K[x_1,\ldots,x_n]_{\bf m}/(\partial f/\partial
x_1,\ldots,\partial f/\partial x_n)=0, 
\end{equation}
hence
\[ \dim_K K[x_1,\ldots,x_n]_{\bf m}/(\partial f/\partial
x_1,\ldots,\partial f/\partial x_n)<\infty. \]
This dimension is called the {\it Milnor number\/} $\mu$ of the
isolated singularity.  We note that (2.1) implies that $\partial
f/\partial x_1,\ldots,\partial f/\partial x_n$ generate an ${\bf
  m}$-primary ideal in $K[x_1,\ldots,x_n]_{\bf m}$, hence form a
regular sequence in that ring.

When ${\rm char}\;K=0$, this definition of isolated singularity is
equivalent to the condition that the hypersurface $f=0$ be nonsingular
in a punctured Zariski neighborhood of ${\bf 0}$ on that hypersurface,
i.~e., that 
\begin{equation}
\mbox{Krull dim}\;K[x_1,\ldots,x_n]_{\bf m}/(f,\partial f/\partial
x_1,\ldots,\partial f/\partial x_n)=0.
\end{equation}
For by the Theorem of Sard-Bertini, the hypersurface $f=c$ is
nonsingular except for finitely many $c\in K$, hence by omitting
finitely many hypersurfaces one obtains a Zariski neighborhood of
${\bf 0}$ in $K^n$ in which ${\bf 0}$ is the only critical point of
the map~$f$.  If ${\rm char}\;K=p>0$, the Theorem of Sard-Bertini
fails and condition (2.2) does not imply that ${\bf 0}$ is an isolated
singularity.  For example, take $f=x_1^p+x_2^a$ with $(p,a)=1$.
Then (2.2) holds but $f$ has infinitely many critical points in any
Zariski neighborhood of ${\bf 0}$ in $K^2$ and 
\[ \mbox{Krull dim}\;K[x_1,x_2]_{\bf m}/(\partial f/\partial
x_1,\partial f/\partial x_2)=1. \]

Recall that $g\in K[x_1,\ldots,x_n]$ is called {\it weighted
  homogeneous of total degree $\delta$} if there exist positive
integers $\alpha_1,\ldots,\alpha_n$ with greatest common divisor $1$
such that
\[ g(\lambda^{\alpha_1}x_1,\ldots,\lambda^{\alpha_n}x_n) =
  \lambda^\delta g(x_1,\ldots,x_n). \]
When this holds, we also have the Euler-type relation
\begin{equation}
\delta g=\sum_{i=1}^n \alpha_ix_i\frac{\partial g}{\partial
  x_i}. 
\end{equation}
Note that $\delta$ and the $\alpha_i$ may not be uniquely determined. 
For example, if $g(x_1,x_2)=x_1x_2$, then $g(\lambda
x_1,\lambda x_2)=\lambda^2 g(x_1,x_2)$ and $g(\lambda x_1,\lambda^2
x_2)=\lambda^3 g(x_1,x_2)$.  

We say that the isolated singularity ${\bf 0}$ of the hypersurface
$f=0$ is {\it weighted homogeneous\/} if there exists a weighted
homogeneous polynomial $g$ such that
\[ K[[x_1,\ldots,x_n]]/(f)\simeq K[[x_1,\ldots,x_n]]/(g). \]
A total degree $\delta$ of $g$ is called a {\it total degree\/} of
the isolated singularity ${\bf 0}$.  In this situation, there exists a
regular system of parameters $x'_1,\ldots,x'_n\in K[[x_1,\ldots,x_n]]$
(i.~e., $n$ elements of $K[[x_1,\ldots,x_n]]$ that generate its
maximal ideal) such that
\begin{equation}
f(x_1,\ldots,x_n)=g(x'_1,\ldots,x'_n). 
\end{equation}
This follows from \cite[Lemma 1.7]{L}, whose proof is valid over an
arbitrary field.

Note that (2.1) implies that when ${\bf 0}$ is an isolated
singularity of $f=0$, there exists a positive integer $m$ such that
\[ f^m\in(\partial f/\partial x_1,\ldots,\partial f/\partial x_n) \]
in the local ring $K[x_1,\ldots,x_n]_{\bf m}$.

\begin{lemma}
Suppose ${\bf 0}$ is a weighted homogeneous isolated singularity of
the hypersurface $f=0$.  If ${\rm char}\;K=p>0$, assume also that
$(p,\delta)=1$, where $\delta$ is a total degree of ${\bf 0}$.  Then
\[ f\in(\partial f/\partial x_1,\ldots,\partial f/\partial x_n) \]
in the local ring $K[x_1,\ldots,x_n]_{\bf m}$.  Furthermore, in every
representation
\[ f=\sum_{i=1}^n h_i\frac{\partial f}{\partial x_i} \]
in $K[x_1,\ldots,x_n]_{\bf m}$, $h_1,\ldots,h_n$ must lie in the
maximal ideal of $K[x_1,\ldots,x_n]_{\bf m}$.
\end{lemma}

{\it Proof}.  From (2.3) it follows that
\[ g\in(\partial g/\partial x_1,\ldots,\partial g/\partial x_n) \]
in $K[x_1,\ldots,x_n]$.  Equation (2.4) then implies that 
\[ f\in (\partial f/\partial x_1,\ldots,\partial f/\partial x_n) \]
in $K[[x_1,\ldots,x_n]]$.  But the natural inclusion
\[ K[x_1,\ldots,x_n]_{\bf m}\hookrightarrow K[[x_1,\ldots,x_n]] \] 
induces an isomorphism
\[ K[x_1,\ldots,x_n]_{\bf m}/(\partial f/\partial x_1,\ldots,\partial
f/\partial x_n)\simeq K[[x_1,\ldots,x_n]]/(\partial f/\partial
x_1,\ldots,\partial f/\partial x_n).  \]
This implies the first assertion of the lemma.  Suppose we have a
representation
\[ f=\sum_{i=1}^n h_i\frac{\partial f}{\partial x_i} \]
in $K[x_1,\ldots,x_n]_{\bf m}$.  By (2.3) and (2.4), we know there is
a representation
\[ f=\sum_{i=1}^n \tilde{h}_i\frac{\partial f}{\partial x_i} \]
with $\tilde{h}_1,\ldots,\tilde{h}_n$ lying in the maximal ideal of
$K[[x_1,\ldots,x_n]]$.  It follows that
\[ \sum_{i=1}^n \frac{\partial f}{\partial x_i}(h_i-\tilde{h}_i)=0. \]
But since $\partial f/\partial x_1,\ldots,\partial f/\partial x_n$
form a regular sequence in $K[x_1,\ldots,x_n]_{\bf m}$, they also form
a regular sequence in $K[[x_1,\ldots,x_n]]$.  Thus there exists a
skew-symmetric set $\{\eta_{ij}\}_{i,j=1}^n\subset
K[[x_1,\ldots,x_n]]$ (i.~e., $\eta_{ji}=-\eta_{ij}$) such that
\[ h_i-\tilde{h}_i=\sum_{j=1}^n \frac{\partial f}{\partial
  x_j}\eta_{ij}. \]
But this implies that $h_1,\ldots,h_n$ lie in the maximal ideal of
$K[[x_1,\ldots,x_n]]$, hence they must also lie in the maximal ideal
of $K[x_1,\ldots,x_n]_{\bf m}$.
\section{Some reduction steps}

We begin with some general remarks on the spectral sequence
$E^{r,s}_t$.  Let
$\omega\in\Omega^m_{{\bf F}_q[x]/{\bf F}_q}$ for some $m$, $0\leq
m\leq n-1$.  If $\omega\in F_r\Omega^m_{{\bf F}_q[x]/{\bf F}_q}$, we
may write
\[ \omega=\sum_{k=0}^r\omega^{(k)}, \]
where $\omega^{(k)}$ is a homogeneous form of degree $k$.  The
assertion that $E^{r,m-r}_e=0$ means that if
\begin{equation}
\sum_{j=0}^i df^{(\delta-j)}\wedge\omega^{(r-i+j)}=0
\end{equation}
for $i=0,1,\ldots,e-1$, then there exist
$\{\xi_j^{(r-\delta+j)}\}_{j=1}^e\subseteq \Omega^{m-1}_{{\bf
    F}_q[x]/{\bf F}_q}$, where $\xi_j^{(r-\delta+j)}$ is homogeneous of
degree $r-\delta+j$, such that
\begin{equation}
\omega^{(r)}=\sum_{j=0}^{e-1} df^{(\delta-j)}\wedge\xi_{j+1}^{(r-\delta+j+1)}
\end{equation}
and such that
\begin{equation}
\sum_{j=0}^i df^{(\delta-j)}\wedge\xi_{j+e-i}^{(r-\delta+j+e-i)}=0
\end{equation}
for $i=0,1,\ldots,e-2$.

Now fix $e=\delta-\delta'+1$.  From (1.7), $f^{(\delta-j)}=0$ for
$0<j<e-1$, thus (3.1) implies
\begin{equation}
df^{(\delta)}\wedge\omega^{(r)}=0
\end{equation}
and
\begin{equation}
df^{(\delta)}\wedge\omega^{(r-\delta+\delta')} +
df^{(\delta')}\wedge\omega^{(r)} =0 
\end{equation}
and (3.2) and (3.3) become
\begin{equation}
\omega^{(r)}=df^{(\delta)}\wedge\xi_1^{(r-\delta+1)} +
df^{(\delta')}\wedge\xi_e^{(r-\delta'+1)}
\end{equation}
and
\begin{equation}
df^{(\delta)}\wedge\xi_e^{(r-\delta'+1)}=0.
\end{equation}

The vanishing of $E_{\delta-\delta'+1}^{r,m-r}$ for all $m<n$ and
all $r$ is thus a consequence of the following stronger assertion.
\begin{proposition}
Assume the hypothesis of Theorem $1.10$.  If $\omega\in \Omega^m_{{\bf
    F}_q[x]/{\bf F}_q}$, $m<n$, is a homogeneous form satisfying
\begin{equation}
df^{(\delta)}\wedge\omega=0
\end{equation}
and 
\begin{equation}
df^{(\delta')}\wedge\omega=df^{(\delta)}\wedge\xi
\end{equation}
for some homogeneous form $\xi\in\Omega^m_{{\bf F}_q[x]/{\bf F}_q}$,
then there exists a homogeneous form $\eta\in\Omega^{m-1}_{{\bf
    F}_q[x]/{\bf F}_q}$ such that
\begin{equation}
\omega=df^{(\delta)}\wedge\eta.
\end{equation}
\end{proposition}

{\it Proof}.  The conclusion for $m<n-1$ follows simply from the fact
that $f^{(\delta)}=0$ has only isolated singularities in ${\bf P}^{n-1}$.  
This implies that the ideal of ${\bf F}_q[x_1,\ldots,x_n]$ generated
by $\{\partial f^{(\delta)}/\partial x_i\}_{i=1}^n$ has height $n-1$,
therefore also has depth $n-1$.  It then follows directly from
\cite{S} that condition (3.9) alone implies the existence of the desired
$\eta$ satisfying (3.11).  In other words, we have
\begin{equation}
H^m(\Omega^{\updot}_{{\bf F}_q[x]/{\bf
    F}_q},\phi_{f^{(\delta)}})=0\qquad\text{for $m<n-1$.}
\end{equation}

So assume that $\omega\in\Omega^{n-1}_{{\bf F}_q[x]/{\bf F}_q}$ is a
homogeneous form satisfying (3.9) and (3.10) for some homogeneous form
$\xi\in\Omega^{n-1}_{{\bf F}_q[x]/{\bf F}_q}$.  We express (3.9) and
(3.10) in coordinate form.  Let 
\begin{align*}
\omega &=\sum_{i=1}^n (-1)^{i-1}\omega_i\,dx_1\wedge\cdots\wedge
\widehat{dx}_i\wedge\cdots\wedge dx_n, \\
\xi &=\sum_{i=1}^n (-1)^{i-1}\xi_i\,dx_1\wedge\cdots\wedge
\widehat{dx}_i\wedge\cdots\wedge dx_n,
\end{align*}
where $\omega_i,\xi_i\in {\bf F}_q[x_1,\ldots,x_n]$ are homogeneous
polynomials.  Then (3.9) becomes
\begin{equation}
\sum_{i=1}^n \frac{\partial f^{(\delta)}}{\partial x_i}\omega_i=0
\end{equation}
and (3.10) becomes
\begin{equation}
\sum_{i=1}^n \frac{\partial f^{(\delta')}}{\partial x_i}\omega_i =
\sum_{i=1}^n \frac{\partial f^{(\delta)}}{\partial x_i}\xi_i. 
\end{equation}

To simplify the calculation, we make a coordinate change.  Let
$a_1,\ldots,a_s\in{\bf P}^{n-1}$ be the singular points of
$f^{(\delta)}=0$.  Since the generic hyperplane section of a
hypersurface with isolated singularities is smooth, we can make a
coordinate change on ${\bf A}^n$ so that the hyperplane $x_n=0$ in
${\bf P}^{n-1}$ intersects the hypersurface $f^{(\delta)}=0$ in ${\bf
  P}^{n-1}$ transversally, in particular, the singularities
$a_1,\ldots,a_s$ do not lie on $x_n=0$.  This implies that the
polynomials $x_n,\partial f^{(\delta)}/\partial x_1,\ldots,\partial
f^{(\delta)}/\partial x_{n-1}$ taken in any order form a regular
sequence in ${\bf F}_q[x_1,\ldots,x_n]$. (We are using here the
hypothesis that $(p,\delta)=1$.) 

We claim that it is enough to show that
\begin{equation}
\omega_n\in\biggl(\frac{\partial f^{(\delta)}}{\partial
  x_1},\ldots,\frac{\partial f^{(\delta)}}{\partial x_{n-1}}\biggr) 
\end{equation}
in ${\bf F}_q[x_1,\ldots,x_n]$.  To see this, suppose
\[ \omega_n=\sum_{i=1}^{n-1} h_i\frac{\partial f^{(\delta)}}{\partial
  x_i} \]
for some homogeneous polynomials $h_i$ and substitute into (3.13) to get
\[ \sum_{i=1}^{n-1} \frac{\partial f^{(\delta)}}{\partial x_i}
\biggl(\omega_i+h_i\frac{\partial f^{(\delta)}}{\partial
  x_n}\biggr)=0. \]
Since $\partial f^{(\delta)}/\partial x_1,\ldots,\partial
  f^{(\delta)}/\partial x_{n-1}$ form a regular sequence, there exists
  a skew-symmetric set $\{\eta_{ij}\}_{i,j=1}^{n-1}$ of homogeneous
  polynomials such that 
\[ \omega_i+h_i\frac{\partial f^{(\delta)}}{\partial
  x_n}=\sum_{j=1}^{n-1} \eta_{ij}\frac{\partial f^{(\delta)}}{\partial x_j}
  \qquad \mbox{for $i=1,\ldots,n-1$.} \]
If we set $\eta_{in}=-h_i$, $\eta_{ni}=h_i$, for $i=1,\ldots,n-1$ and
  $\eta_{nn}=0$, then $\{\eta_{ij}\}_{i,j=1}^n$ is a skew-symmetric
  set satisfying
\begin{equation}
\omega_i=\sum_{j=1}^n \eta_{ij}\frac{\partial f^{(\delta)}}{\partial
  x_j} \qquad \mbox{for $i=1,\ldots,n$.}
\end{equation}
If we then define
\[ \eta=\sum_{1\leq i<j\leq n} (-1)^i\eta_{ij}\,
dx_1\wedge\cdots\wedge\widehat{dx}_i \wedge\cdots\wedge\widehat{dx}_j
\wedge\cdots\wedge dx_n, \]
equation (3.16) implies equation (3.11).

The common zeros of $\partial f^{(\delta)}/\partial x_1,\ldots,\partial
f^{(\delta)}/\partial x_{n-1}$ in ${\bf P}^{n-1}$ form a finite set
containing the singular points of $f^{(\delta)}=0$, so we may write this
set as
\[ \{a_1,\ldots,a_s,b_1,\ldots,b_t\}. \]
Since $x_n,\partial f^{(\delta)}/\partial x_1,\ldots,\partial
f^{(\delta)}/\partial x_{n-1}$ form a regular sequence, none of these
points lies on the hypersurface $x_n=0$.  Note that our hypotheses
imply that the hypersurface $\partial f^{(\delta)}/\partial x_n=0$ in ${\bf
  P}^{n-1}$ contains the points $a_1,\ldots,a_s$ but does not contain
any of the points $b_1,\ldots,b_t$.  The main technical tool for
proving (3.15) is the following.

\begin{lemma}
There exists a homogeneous polynomial $P\in {\bf F}_q[x_1,\ldots,x_n]$
such that
\[ P\omega_n\in\biggl(\frac{\partial f^{(\delta)}}{\partial
  x_1},\ldots,\frac{\partial f^{(\delta)}}{\partial x_{n-1}}\biggr) \]
and such that the hypersurface $P=0$ in ${\bf P}^{n-1}$ does not
contain any of the points $a_1,\ldots,a_s$.  
\end{lemma}

{\it Remark}.  The proof of Lemma 3.17 will require several steps.
Before starting the proof, we explain how it implies (3.15).  For any fixed
$j\in\{1,\ldots,t\}$, we can find a linear form $h_j\in{\bf
  F}_q[x_1,\ldots,x_n]$ such that the hyperplane $h_j=0$ in ${\bf
  P}^{n-1}$ contains $b_j$ but contains none of $a_1,\ldots,a_s$.
Multiplying $P$ by such factors, we may assume in addition to the
conclusion of the lemma that the hypersurface $P=0$ in ${\bf P}^{n-1}$
contains $b_1,\ldots,b_t$.  Choose nonnegative integers
$\alpha,\beta$ such that $x_n^{\alpha}P+x_n^{\beta}\partial
f^{(\delta)}/\partial x_n$ is homogeneous.  The properties of $P$
imply that the hypersurface $x_n^{\alpha}P+x_n^{\beta}\partial
f^{(\delta)}/\partial x_n=0$ in ${\bf P}^{n-1}$ contains none of the points
$a_1,\ldots,a_s,b_1,\ldots,b_t$.  Thus the homogeneous polynomials 
\[ \frac{\partial f^{(\delta)}}{\partial x_1},\ldots,\frac{\partial
  f^{(\delta)}}{\partial x_{n-1}}, x_n^{\alpha}P+x_n^{\beta}\frac{\partial
  f^{(\delta)}}{\partial x_n} \] 
have no common zero in ${\bf P}^{n-1}$ and hence form a regular
sequence in ${\bf F}_q[x_1,\ldots,x_n]$.  But (3.13) and Lemma 3.17
  imply that 
\[ \biggl(x_n^{\alpha}P+x_n^{\beta}\frac{\partial f^{(\delta)}}{\partial
  x_n}\biggr)\omega_n\in\biggl(\frac{\partial f^{(\delta)}}{\partial 
  x_1},\ldots,\frac{\partial f^{(\delta)}}{\partial x_{n-1}}\biggr). \]
This implies (3.15).
\section{Proof of Lemma 3.17}

There are two basic ideas involved in the proof of Lemma 3.17.  The
first is expressed in the following.

\begin{lemma}
For each singular point $a_i$, $i=1,\ldots,s$, there exist homogeneous
polynomials $Q_i,R^{(i)}_1,\ldots,R^{(i)}_{n-1}$ such that
\begin{equation}
Q_if^{(\delta)}=\sum_{j=1}^{n-1} R^{(i)}_j\frac{\partial
  f^{(\delta)}}{\partial x_j} 
\end{equation}
and such that $a_i$ does not lie on the hypersurface $Q_i=0$ in ${\bf
  P}^{n-1}$ but does lie on all the hypersurfaces $R^{(i)}_j=0$ for
$j=1,\ldots,n-1$.  
\end{lemma}

{\it Proof}.  Fix $i$ and let $(\tilde{a}_1,\ldots,\tilde{a}_{n-1},1)$
be homogeneous coordinates for $a_i\in{\bf P}^{n-1}$.  Put
\[ \tilde{f}(y_1,\ldots,y_{n-1})=f^{(\delta)}(y_1,\ldots,y_{n-1},1). \]
By Lemma 2.5 (we are using here the hypothesis that
$(p,\delta_i)=1$), we have 
\begin{equation}
\tilde{f}=\sum_{j=1}^{n-1}\tilde{h}_j\frac{\partial
  \tilde{f}}{\partial y_j}, 
\end{equation}
where $\tilde{h}_1,\ldots,\tilde{h}_{n-1}$ lie in the maximal ideal of
the local ring of $(\tilde{a}_1,\ldots,\tilde{a}_{n-1})$, i.~e.,
$\tilde{h}_j=\tilde{P}_j/\tilde{Q}_j$ where
$\tilde{P}_j,\tilde{Q}_j\in K[y_1,\ldots,y_{n-1}]$ and 
\begin{eqnarray*}
\tilde{Q}_j(\tilde{a}_1,\ldots,\tilde{a}_{n-1}) & \neq & 0 \\
\tilde{P}_j(\tilde{a}_1,\ldots,\tilde{a}_{n-1}) & = & 0
\end{eqnarray*}
for $j=1,\ldots,n-1$.  Multiplying (4.3) by
$\tilde{Q}:=\tilde{Q}_1\cdots\tilde{Q}_{n-1}$ gives a relation 
\begin{equation}
\tilde{Q}\tilde{f}=\sum_{j=1}^{n-1} \tilde{R}_j\frac{\partial
  \tilde{f}}{\partial y_j} 
\end{equation}
in ${\bf F}_q[y_1,\ldots,y_{n-1}]$ with
$\tilde{Q}(\tilde{a}_1,\ldots,\tilde{a}_{n-1})\neq 0$ and
$\tilde{R}_j(\tilde{a}_1,\ldots,\tilde{a}_{n-1})=0$ for 
$j=1,\ldots,n-1$.  Making the substitution $y_j\mapsto x_j/x_n$ in
(4.4) and multiplying by a sufficiently high power of $x_n$ then gives
the desired assertion. 

By the argument used in the remark following Lemma 3.17, we may assume
in addition to the conclusion of Lemma 4.1 that the hypersurfaces
$Q_i=0$, $R^{(i)}_1=0$,\ldots, $R^{(i)}_{n-1}=0$ in ${\bf P}^{n-1}$
all contain the points $a_1,\ldots,\hat{a}_i,\ldots,a_s$.
Choose nonnegative integers $\alpha_1,\ldots,\alpha_s$ such that
\[ Q=x_n^{\alpha_1}Q_1+\cdots+x_n^{\alpha_s}Q_s \]
is homogeneous.  Multiplying (4.2) by $x_n^{\alpha_i}$ and summing over
$i$ then gives the following.

\begin{corollary}
There exist homogeneous polynomials $Q,R_1,\ldots,R_{n-1}$ such that
\[ Qf^{(\delta)}=\sum_{j=1}^{n-1} R_j\frac{\partial
  f^{(\delta)}}{\partial x_j} \]
and such that the hypersurface $Q=0$ in ${\bf P}^{n-1}$ contains none
of the points $a_1,\ldots,a_s$ and the hypersurfaces $R_j=0$ contain
the points $a_1,\ldots,a_s$ for $j=1,\ldots,n-1$.
\end{corollary}

Multiplying the Euler relation for $f^{(\delta)}$ by $Q$ gives
\[ \delta Qf^{(\delta)}=\sum_{j=1}^n x_jQ\frac{\partial
  f^{(\delta)}}{\partial x_j}. \]
Combined with Corollary 4.5, this gives
\begin{equation}
x_nQ\frac{\partial f^{(\delta)}}{\partial x_n}=\sum_{j=1}^{n-1}
S_j\frac{\partial f^{(\delta)}}{\partial x_j},
\end{equation}
where
\begin{equation}
S_j=\delta R_j-x_jQ \qquad \mbox{for $j=1,\ldots,n-1$.}
\end{equation}

We can now prove Lemma 3.17.  Multiplying (3.13) by $x_nQ$ and using
(4.6) leads to
\[ \sum_{j=1}^{n-1} \frac{\partial f^{(\delta)}}{\partial
  x_j}(x_nQ\omega_j+S_j\omega_n)=0. \]
Since $\partial f^{(\delta)}/\partial x_1,\ldots,\partial
f^{(\delta)}/\partial x_{n-1}$ form a regular sequence, there exists a
skew-symmetric set $\{\eta_{ij}\}_{i,j=1}^{n-1}$ of homogeneous
polynomials such that 
\begin{equation}
x_nQ\omega_j+S_j\omega_n=\sum_{k=1}^{n-1}\eta_{jk}
\frac{\partial f^{(\delta)}}{\partial x_k} \qquad \mbox{for
  $j=1,\ldots,n-1$.}
\end{equation}
Multiplying (3.14) by $x_nQ$ and using (4.6) gives
\[ \sum_{j=1}^n \frac{\partial f^{(\delta')}}{\partial x_j}
x_nQ\omega_j\in \biggl(\frac{\partial f^{(\delta)}}{\partial
  x_1},\ldots,\frac{\partial f^{(\delta)}}{\partial x_{n-1}}\biggr). \]
Substitution from (4.8) then gives
\begin{equation}
\biggl(-\sum_{j=1}^{n-1}\frac{\partial f^{(\delta')}}{\partial
  x_j}S_j+\frac{\partial f^{(\delta')}}{\partial
  x_n}x_nQ\biggr)\omega_n\in \biggl(\frac{\partial f^{(\delta)}}{\partial
  x_1},\ldots,\frac{\partial f^{(\delta)}}{\partial x_{n-1}}\biggr).
\end{equation}

We now come to the second basic idea of the proof.  Put
\[ P=-\sum_{j=1}^{n-1}\frac{\partial f^{(\delta')}}{\partial
  x_j}S_j+\frac{\partial f^{(\delta')}}{\partial x_n}x_nQ, \]
a homogeneous polynomial.  By (4.9), $P$ satisfies the first assertion
of Lemma~3.17.  We show that it satisfies the second assertion as
well.  Let $(c_1,\ldots,c_n)$ be a set of homogeneous coordinates for
one of the points $a_1,\ldots,a_s$.  By (4.7) and
Corollary~4.5, we see that
\begin{eqnarray*}
P(c_1,\ldots,c_n) & = & \sum_{j=1}^n
c_jQ(c_1,\ldots,c_n)\frac{\partial f^{(\delta')}}{\partial
  x_j}(c_1,\ldots,c_n) \\
 & = & \delta' Q(c_1,\ldots,c_n)f^{(\delta')}(c_1,\ldots,c_n)
\end{eqnarray*}
using the Euler relation for $f^{(\delta')}$.  By hypothesis
$\delta' f^{(\delta')}(c_1,\ldots,c_n)\neq 0$ and by Corollary 4.5 
$Q(c_1,\ldots,c_n)\neq 0$, hence $P(c_1,\ldots,c_n)\neq 0$.  This
proves Lemma~3.17, which completes the proof of Theorem 1.10. 

\section{Proof of Theorem 1.11}

Throughout this section, we assume the hypothesis of Theorem 1.11.
Since the generic hyperplane section of a hypersurface with isolated
singularities is smooth, we may assume, after a coordinate change if
necessary, that the hyperplane $x_n=0$ intersects $f^{(\delta)}=0$
transversally.  Let $\tilde{f}\in{\bf F}_q[x_1,\ldots,x_{n-1}]$ be
defined by 
\[ \tilde{f}(x_1,\ldots,x_{n-1})=f^{(\delta)}(x_1,\ldots,x_{n-1},0). \]
Then $\tilde{f}=0$ defines a smooth hypersurface in ${\bf P}^{n-2}$.

\begin{lemma}
Under the above conditions,  $\partial f^{(\delta)}/\partial
x_1,\ldots,\partial f^{(\delta)}/\partial x_{n-1}$ form a regular
sequence. 
\end{lemma}

{\it Proof}.  It suffices to show that $\{\partial f^{(\delta)}/\partial
x_i\}_{i=1}^{n-1}$ have only finitely many common zeros in ${\bf
  P}^{n-1}$.  Since $p|\delta$, the Euler relation becomes
\begin{equation}
x_n\frac{\partial f^{(\delta)}}{\partial x_n}=-\sum_{i=1}^{n-1}
x_i\frac{\partial f^{(\delta)}}{\partial x_i}, 
\end{equation}
thus any common zero of $\{\partial f^{(\delta)}/\partial
x_i\}_{i=1}^{n-1}$ is a zero of either $\partial f^{(\delta)}/\partial x_n$
or $x_n$.  Those which are zeros of $\partial f^{(\delta)}/\partial x_n$
form a finite set by the hypothesis of Theorem 1.11.  Those which are
zeros of $x_n$ are in one-to-one correspondence with the zeros of
$\{\partial\tilde{f}/\partial x_i\}_{i=1}^{n-1}$ in~${\bf P}^{n-2}$.
Since $\tilde{f}=0$ defines a smooth hypersurface in ${\bf P}^{n-2}$,
this set must also be finite.

To prove Theorem 1.11, the discussion in section 3 shows that it
suffices to prove the analogue of Proposition 3.8.
\begin{proposition}
Assume the hypothesis of Theorem $1.11$.  If $\omega\in \Omega^m_{{\bf
    F}_q[x]/{\bf F}_q}$, $m<n$, is a homogeneous form satisfying
\begin{equation}
df^{(\delta)}\wedge\omega=0
\end{equation}
and 
\begin{equation}
df^{(\delta')}\wedge\omega=df^{(\delta)}\wedge\xi
\end{equation}
for some homogeneous form $\xi\in\Omega^m_{{\bf F}_q[x]/{\bf F}_q}$,
then there exists a homogeneous form $\eta\in\Omega^{m-1}_{{\bf
    F}_q[x]/{\bf F}_q}$ such that
\begin{equation}
\omega=df^{(\delta)}\wedge\eta.
\end{equation}
\end{proposition}

{\it Proof}.  By Lemma 5.1, the ideal of ${\bf F}_q[x_1,\ldots,x_n]$
generated by $\{\partial f^{(\delta)}/\partial x_i\}_{i=1}^n$ has
depth $n-1$.  Thus, just as in the proof of Proposition 3.8, we
conclude that for $m<n-1$, (5.4) alone implies (5.6).  Hence (3.12)
holds under the hypothesis of Theorem 1.11 also.

Suppose now $m=n-1$.  As in the proof of Proposition 3.8, we write
(5.4) and (5.5) in coordinate form:
\begin{align}
\sum_{i=1}^n \frac{\partial f^{(\delta)}}{\partial x_i}\omega_i &= 0 \\
\sum_{i=1}^n \frac{\partial f^{(\delta')}}{\partial x_i}\omega_i &=
\sum_{i=1}^n \frac{\partial f^{(\delta)}}{\partial x_i}\xi_i.
\end{align}
The same argument as before (see (3.15)) reduces us to proving that 
\begin{equation}
\omega_n\in\biggl(\frac{\partial f^{(\delta)}}{\partial
  x_1},\ldots,\frac{\partial f^{(\delta)}}{\partial x_{n-1}}\biggr).
\end{equation}

\begin{lemma}
Suppose that $(5.7)$ holds and that
\[ f^{(\delta')}\omega_n\in\biggl(\frac{\partial f^{(\delta)}}{\partial
x_1},\ldots,\frac{\partial f^{(\delta)}}{\partial x_{n-1}}\biggr). \]
Then
\[ \omega_n\in\biggl(\frac{\partial f^{(\delta)}}{\partial
  x_1},\ldots,\frac{\partial f^{(\delta)}}{\partial x_{n-1}}\biggr). \]
\end{lemma}

{\it Proof}.  The zeros of $\{\partial f^{(\delta)}/\partial
x_i\}_{i=1}^n$ in ${\bf P}^{n-1}$ form a finite set
$\{a_1,\ldots,a_s\}$ and the zeros of $\{\partial
f^{(\delta)}/\partial x_i\}_{i=1}^{n-1}$ form a finite set
$\{a_1,\ldots,a_s,b_1,\ldots,b_t\}$.  For $i=1,\ldots,t$, choose a
linear form $h_i\in{\bf F}_q[x_1,\ldots,x_n]$ that vanishes at $b_i$
but not at $a_j$ for any $j$.  Let $k\in{\bf F}_q[x_1,\ldots,x_n]$ be
a linear form that does not vanish at any $b_i$.  For suitably chosen
nonnegative integers $\alpha$ and $\beta$, the polynomial
\[ g:=(h_1\cdots h_t)^{\alpha}f^{(\delta')}+k^{\beta}\frac{\partial
  f^{(\delta)}}{\partial x_n}\in{\bf F}_q[x_1,\ldots,x_n] \]
is homogeneous.  By the hypothesis of Theorem 1.11, $f^{(\delta')}$
does not vanish at $a_i$ for any $i$.  It follows that $g$ does not
vanish at any $a_i$ or $b_j$, i.~e., the polynomials
\[ \frac{\partial f^{(\delta)}}{\partial x_1},\ldots,\frac{\partial
  f^{(\delta)}}{\partial x_{n-1}},g \]
have no common zero in ${\bf P}^{n-1}$, hence they form a regular
sequence.  But the hypothesis of the lemma implies that
\[ g\omega_n\in(\partial f^{(\delta)}/\partial x_1,\ldots,\partial
f^{(\delta)}/\partial x_{n-1}). \]
The conclusion of the lemma now follows from the defining property of
regular sequences.

By Lemma 5.10, we are reduced to showing the following.
\begin{lemma}
If $(5.7)$ and $(5.8)$ hold, then
\[ f^{(\delta')}\omega_n\in \biggl(\frac{\partial f^{(\delta)}}{\partial
  x_1},\ldots,\frac{\partial f^{(\delta)}}{\partial x_{n-1}}\biggr). \]
\end{lemma}

{\it Proof}.  Multiplying (5.7) by $x_n$ and substituting from the Euler
relation (5.2) gives
\[ \sum_{i=1}^{n-1} \frac{\partial f^{(\delta)}}{\partial
  x_i}(x_n\omega_i - x_i\omega_n)=0. \]
By Lemma 5.1, this implies
\begin{equation}
x_n\omega_i-x_i\omega_n\in \biggl(\frac{\partial
  f^{(\delta)}}{\partial x_1},\ldots,\frac{\partial f^{(\delta)}}{\partial
  x_{n-1}}\biggr) 
\end{equation}
for $i=1,\ldots,n$.  Multiplying (5.8) by $x_n$ and using (5.2)
  gives
\begin{equation}
\sum_{i=1}^n x_n\omega_i\frac{\partial f^{(\delta')}}{\partial
  x_i}\in \biggl(\frac{\partial f^{(\delta)}}{\partial
  x_1},\ldots,\frac{\partial f^{(\delta)}}{\partial x_{n-1}}\biggr).
\end{equation}
It follows from (5.12) and (5.13) that
\[ \sum_{i=1}^n x_i\omega_n\frac{\partial f^{(\delta')}}{\partial
  x_i}\in \biggl(\frac{\partial f^{(\delta)}}{\partial
  x_1},\ldots,\frac{\partial f^{(\delta)}}{\partial x_{n-1}}\biggr). \]
The Euler relation for $f^{(\delta')}$ now implies
\[ \delta' f^{(\delta')}\omega_n\in\biggl(\frac{\partial
  f^{(\delta)}}{\partial x_1},\ldots,\frac{\partial
  f^{(\delta)}}{\partial x_{n-1}}\biggr). \] 
The conclusion of the lemma then follows from the hypothesis that
$(p,\delta')=1$.

\section{Formula for $M_f$}

By \cite[section 1]{AS3}, we know that if
$E^{r,s}_{\delta-\delta'+1}=0$ for all $r,s$ with $r+s\neq n$, then
\begin{equation}
M_f=\dim_{{\bf F}_q}\biggl(\bigoplus_{r+s=n}
E^{r,s}_{\delta-\delta'+1}\biggr).  
\end{equation}
We describe the terms on the right-hand side explicitly.

For $0\leq m\leq n$, let $H^m(\Omega^{\updot}_{{\bf F}_q[x]/{\bf
    F}_q},\phi_{f^{(\delta)}})^{(r)}$ denote the homogeneous component of
degree $r$ of $H^m(\Omega^{\updot}_{{\bf F}_q[x]/{\bf
    F}_q},\phi_{f^{(\delta)}})$ relative to the grading on $\Omega^m_{{\bf
    F}_q[x]/{\bf F}_q}$ defined in section 1.  We define a map
\begin{equation}
\phi_{f^{(\delta')}}:H^m(\Omega^{\updot}_{{\bf F}_q[x]/{\bf
    F}_q},\phi_{f^{(\delta)}})\rightarrow H^{m+1}(\Omega^{\updot}_{{\bf
    F}_q[x]/{\bf F}_q},\phi_{f^{(\delta)}}).
\end{equation}
Let $\omega\in\Omega^m_{{\bf F}_q[x]/{\bf F}_q}$ be such
that $df^{(\delta)}\wedge\omega=0$ and let $[\omega]\in
H^m(\Omega^{\updot}_{{\bf F}_q[x]/{\bf F}_q},\phi_{f^{(\delta)}})$ be the
cohomology class of $\omega$.  We define
\begin{equation}
\phi_{f^{(\delta')}}([\omega])=[df^{(\delta')}\wedge\omega].
\end{equation}
From the definition of the spectral sequence $E^{r,s}_t$, one sees
that $E^{r,n-r}_{\delta-\delta'+1}$ is the cokernel of the map
\[ \phi_{f^{(\delta')}}:H^{n-1}(\Omega^{\updot}_{{\bf F}_q[x]/{\bf
    F}_q},\phi_{f^{(\delta)}})^{(r-\delta'+\delta)}\rightarrow
    H^n(\Omega^{\updot}_{{\bf F}_q[x]/{\bf
    F}_q},\phi_{f^{(\delta)}})^{(r)}. \]
It follows from (6.1) that
\begin{equation}
M_f = \dim_{{\bf F}_q}({\rm
    coker}(\phi_{f^{(\delta')}}:H^{n-1}(\Omega^{\updot}_{{\bf F}_q[x]/{\bf 
    F}_q},\phi_{f^{(\delta)}})\rightarrow H^n(\Omega^{\updot}_{{\bf
    F}_q[x]/{\bf F}_q},\phi_{f^{(\delta)}}))).
\end{equation}
We compute the dimension of this cokernel.

Under the hypothesis of either Theorem 1.10 or 1.11, $\{\partial
f^{(\delta)}/\partial x_i\}_{i=1}^n$ have finitely many common zeroes in
${\bf P}^{n-1}$, say, $a_1,\ldots,a_s$.  By a coordinate change, we
may assume $a_1,\ldots,a_s$ lie in the open set $x_n\neq 
0$, which we identify with ${\bf A}^{n-1}$.  Put
\[ h=f^{(\delta)}(y_1,\ldots,y_{n-1},1)\in{\bf F}_q[y_1,\ldots,y_{n-1}]. \]
The Milnor number $\mu_i$ of $a_i$ is given by
\[ \mu_i=\dim_{{\bf F}_q}{\bf F}_q[y_1,\ldots,y_{n-1}]_{{\bf m}_i} /
(\partial h/\partial y_1,\ldots,\partial h/\partial y_{n-1}), \]
where ${\bf F}_q[y_1,\ldots,y_{n-1}]_{{\bf m}_i}$ denotes the
localization of ${\bf F}_q[y_1,\ldots,y_{n-1}]$ at the maximal ideal
${\bf m}_i$ corresponding to $a_i$.

\begin{proposition}
Under the hypothesis of either Theorem $1.10$ or $1.11$, 
\[ M_f=(\delta-1)^n-(\delta-\delta')\sum_{i=1}^s \mu_i. \]
\end{proposition}

{\it Proof}.  The graded module $H^i(\Omega^{\updot}_{{\bf
    F}_q[x]/{\bf F}_q},\phi_{f^{(\delta)}})$ has a Poincar\'{e} series
    $p_i(t)$: 
\[ p_i(t)=\sum_{r=0}^{\infty} \biggl(\dim_{{\bf F}_q}
H^i(\Omega^{\updot}_{{\bf F}_q[x]/{\bf F}_q},\phi_{f^{(\delta)}})^{(r)}
\biggr)t^r. \] 
Using only $\deg f=\delta$, one has always
\[ \sum_{i=0}^n (-1)^{n-i}p_i(t)=\frac{(1-t^{\delta-1})^n}{(1-t)^n}. \]
Under the hypothesis of either Theorem 1.10 or 1.11, equation
(3.12) holds.  Thus $p_i(t)=0$ for $i<n-1$ and we have
\begin{equation}
p_n(t)-p_{n-1}(t)=\frac{(1-t^{\delta-1})^n}{(1-t)^n}.
\end{equation}

Put $f_n=\partial f^{(\delta)}/\partial x_n$ and define
\[ h_n=f_n(y_1,\ldots,y_{n-1},1)\in {\bf F}_q[y_1,\ldots,y_{n-1}]. \]
The Euler relation for $f^{(\delta)}$ implies
\begin{equation}
\delta h=h_n+\sum_{i=1}^{n-1} y_i\frac{\partial h}{\partial y_i}.
\end{equation}
The proof of Choudary-Dimca\cite[Corollary 9]{CD} shows that for all
sufficiently large $r$, 
$\dim_{{\bf F}_q} H^n(\Omega^{\updot}_{{\bf F}_q[x]/{\bf
    F}_q},\phi_{f^{(\delta)}})^{(r)}$ 
is a constant equal to 
\begin{equation}
\sum_{i=1}^s \dim_{{\bf F}_q} {\bf F}_q[y_1,\ldots,y_{n-1}]_{{\bf m}_i}
  / (\partial h/\partial y_1,\ldots,\partial h/\partial y_{n-1},h_n).
\end{equation}
Since the right-hand side of (6.6) is a polynomial in $t$, it follows
that for sufficiently large $r$, $\dim_{{\bf F}_q} 
H^{n-1}(\Omega^{\updot}_{{\bf F}_q[x]/{\bf
    F}_q},\phi_{f^{(\delta)}})^{(r)}$ also equals (6.8).  If
$p|\delta$, (6.7) shows that $h_n$ lies in the ideal generated by
$\{\partial h/\partial y_i\}_{i=1}^{n-1}$.  If $a_i$ is a weighted
homogeneous isolated singular point, then $h$ lies in the ideal of
${\bf F}_q[y_1,\ldots,y_{n-1}]_{{\bf m}_i}$ generated by $\{\partial
h/\partial y_i\}_{i=1}^{n-1}$, so by (6.7) $h_n$ also lies in this
ideal.  Thus in either case, (6.8) equals $\sum_{i=1}^s \mu_i$.  We
can summarize these facts by saying that there exist polynomials
$q_n(t),q_{n-1}(t)$ such that
\begin{equation}
p_n(t)=\frac{q_n(t)}{1-t}, \qquad p_{n-1}(t)=\frac{q_{n-1}(t)}{1-t},
\end{equation}
and such that
\begin{equation}
q_n(1)=q_{n-1}(1)=\sum_{i=1}^s \mu_i.
\end{equation}

By Propositions 3.8 and 5.3 the mapping (6.2) is injective for $m\leq
n-1$, and it is homogeneous of degree $\delta'-\delta$ in the grading we have
defined, hence the Poincar\'{e} series of its cokernel is
\begin{align*}
p_n(t)-t^{\delta'-\delta}p_{n-1}(t) &=
p_n(t)-p_{n-1}(t)+(1-t^{\delta'-\delta})p_{n-1}(t) \\
 &= p_n(t)-p_{n-1}(t)+\frac{t^{\delta-\delta'}-1}{1-t}
\frac{q_{n-1}(t)}{t^{\delta-\delta'}}
\end{align*}
using (6.9).  By (6.6), this expression simplifies to
\[ (1+t+\cdots+t^{\delta-2})^n-(1+t+\cdots+t^{\delta-\delta'-1})
\frac{q_{n-1}(t)}{t^{\delta-\delta'}}. \]
Using (6.10), we see that the
value of this polynomial at $t=1$ is
\[ (\delta-1)^n-(\delta-\delta')\sum_{i=1}^s \mu_i, \]
which completes the proof of Proposition 6.5.

\end{document}